# A Bi-Stage Framework for Automatic Development of Pixel-Based Planar Antenna Structures

Khadijeh Askaripour[1[0000-0003-0487-9502]], Adrian Bekasiewicz[1[0000-0003-0244-541X]], and Slawomir Koziel[1,2[0000-0002-9063-2647]]

[1] Faculty of Electronics, Telecommunications and Informatics, Gdansk University of Technology, Narutowicza 11/12, 80-233 Gdansk, Poland
[2] Department of Engineering, Reykjavik University, Menntavegur 1, 102 Reykjavik, Iceland
`adrian.bekasiewicz@pg.edu.pl`

**Abstract.** Development of modern antennas is a cognitive process that intertwines experience-driven determination of topology and tuning of its parameters to fulfill the performance specifications. Alternatively, the task can be formulated as an optimization problem so as to reduce reliance of geometry selection on engineering insight. In this work, a bi-stage framework for automatic generation of antennas is considered. The method determines free-form topology through optimization of interconnections between components (so-called pixels) that constitute the radiator. Here, the process involves global optimization of connections between pixels followed by fine-tuning of the resulting topology using a surrogate-assisted local-search algorithm to fulfill the design requirements. The approach has been demonstrated based on two case studies concerning development of broadband and dual-band monopole antennas.

## 1 Introduction

Antenna design is a challenging task that involves development of topology followed by its tuning to obtain the desired performance specifications. Given lack of empirical, or analytical solutions that could streamline the process, modern radiators are generated and tuned based on responses obtained from electromagnetic (EM) simulations of their models [1]. EM-driven evolution of topology w.r.t. application-specific requirements is considered a pivotal part of the design process [2]. However, it is also subject to engineering bias which might hinder identification of high performance solutions with unintuitive geometries [1], [3], [4]. From this perspective, automatic design of antennas is an interesting alternative to cognition-driven procedures.

Algorithmic design methods involve EM-driven development of free-form radiators represented using either a set of characteristic points (interconnected using line sections, or splines), or a distribution of primitives (e.g., rectangles) specified in the form of a binary matrix [4]-[7]. However, a large number of parameters required to ensure flexibility of radiator shape makes both representations impractical for automatic antenna design [5], [7]. This is because exploration of large search spaces associated with free-form geometries involves the use of global-search methods (e.g., population-based metaheuristics). The latter ones are numerically prohibitive for multi-dimensional problems, due to excessive number of EM simulations (hundreds or even thousands) required to converge [4]-[7].

For structures represented using primitives, the problem of unacceptable design cost can be mitigated using internal multi-port method (IMPM) [8], [9]. The approach boils down to representation of the radiator as a lattice of dummy components (so-called pixels) interconnected through internal ports [8]. Identification of IMPM model involves extraction of an impedance matrix pertinent to pixel-based antenna from a single, multi-port EM simulation [9]. Electrical performance of the radiator is then approximated through adjustment of impedance between the internal ports. The approach represents a significant advancement compared to conventional models where distribution of primitives (e.g., metallic rectangles) is specified within an associated matrix and structure responses are derived from EM simulations (one at a time) [7]. Instead, IMPM enables adjustment of connections between pixels (and antenna performance characteristics) at a negligible cost during post-processing (simple calculations based on circuit theory). Hence, method can be used to reliably handle global optimization of complex topologies [8]. Regardless of numerical efficiency, IMPM does not support fine-tuning of responses based on adjustment of pixels dimensions.

In this work, a bi-step framework for automatic development of pixel-based antennas is considered. The method involves global optimization of IMPM-based antenna model to approximate the desired solution. The obtained design is then used as a starting point for surrogate-assisted fine-tuning. The approach has been demonstrated using two test cases concerning development of broadband and dual-band antennas based on monopole topology. For the considered numerical experiments, the computational cost of antenna design does not exceed 36 EM model simulations.

## 2      Design Framework

### 2.1      Problem Formulation

Let $R(x, y)$ be the EM simulation response of a generic pixel-based antenna model obtained over a frequency range of interest $f$ for the vector $x$ that represents floating-point dimensions and binary vector $y$ that denotes interconnections between the specific pixels (for visualization, see Fig. 1). The design problem can be defined as the optimization task of the form:

$$x^*, y^* = \arg \min_{x,y} \left( U\left( R(x, y) \right) \right) \qquad (1)$$

where $U$ is a scalar objective function; $x^*$ and $y^*$ represent the optimum design parameters and pixel connections to be found. Direct minimization of (1) is computationally impractical as it would necessitate hundreds, or thousands of $R(x, y)$ simulations to complete. Instead, the design can be realized as a bi-stage process, where appropriate configuration of *y* is first obtained as a result of global-search executed on the IMPM representation of the antenna $R_{sy}(x, y)$ with constant *x* by solving:

$$y^* = \arg\min_{y}\left(U\left(R_{sy}(x, y)\right)\right) \tag{2}$$

Next, the specific dimensions of the radiator *x*—with constant $y^*$ obtained from (2)—can be subject to fine-tuning using a local-search algorithm.

### 2.2   IMPM Representation of Pixel-Based Antenna

Consider an EM-based representation of an $N_x$ row and $N_y$ column pixel-based radiator in monopole configuration $R(x, y)$ and a schematic view of its IMPM model $R_{sy}(y) = R_{sy}(x, y)$, both shown in Fig. 1. The vector of antenna floating-point dimensions is given as $x = [l\ d\ \alpha\ \gamma]^T$. The relative variables are: $l_s = lN_y + d(N_y - 1) + 2o$, $l_m = w_g + \gamma$, $l_g = l_s - \alpha$, and $w_g = \beta l_m$, whereas $w_m = 3$, $\beta = 0.4$, and $o = 3$ remain constant. The unit for all dimensions (except unit-less $\alpha$, $\beta$, $\gamma$) is mm. Parameter $M = 2N_xN_y - (N_x + N_y)$ specifies the number of unique internal (auxiliary) ports between the pixels, whereas $y = [z_{l.2} \ldots z_{l.N}]^T$ is an *M*-element bi-value vector of load impedances for individual ports [9]. Its component associated with *m*th port ($m = 2, \ldots, N$) can be set to either 0, or ∞ when pixels are to be connected (closed port), or disconnected (open port). The IMPM model can be outlined as follows. Let $Z = Z(x, f)$ be an $N \times N$ matrix [8]:

$$Z = \begin{bmatrix} Z_A & Z_B \\ Z_C & Z_D \end{bmatrix} \tag{3}$$

where $Z_A = z_{1.1}$, is pertinent to external port impedance; $Z_B = [z_{1.2} \ldots z_{1.N}]$, $Z_C = [z_{2.1} \ldots z_{N.1}]^T$ represent vectors of external-to-auxiliary ports and $Z_D$ is a matrix of auxiliary ports impedances (cf. Fig. 1(b)). The components of (3) are calculated only once for the given configuration of pixel antenna EM-based model [8]. Note that configuration is understood here as the number of pixels and ports *M* between them, but also as a specific vector of floating point parameters *x* and the frequency sweep *f*. Once (3) is extracted, the impedance at external input port can be obtained by solving [8]:

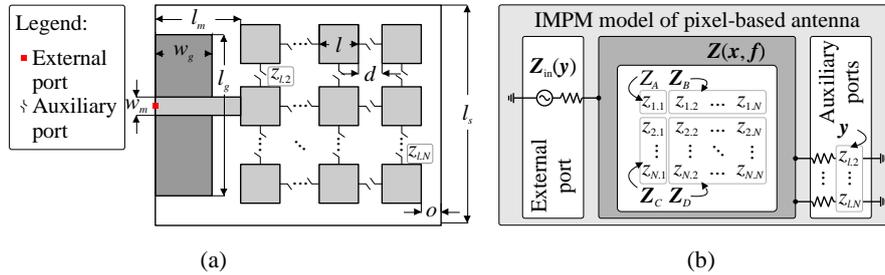

(a)                        (b)

**Fig. 1.** Pixel-based antenna development: (a) visualization of a generic EM model with highlight on dimensions and ports (along with impedances associated with them) and (b) representation of the structure input impedance using IMPM model.

$$Z_{in}(y) = Z_A - Z_B (Y(y) + Z_D)^{-1} Z_C \qquad (4)$$

Here, $Y(y)$ is an $N \times N$ diagonal matrix comprising the elements of $y$ vector. Finally, IMPM model response $R_{sy}(y)$ can be derived from (4) as [8]:

$$R_{sy}(y) = \left(\frac{Z_{in}(y)}{z_0} + I\right)^{-1} \left(\frac{Z_{in}(y)}{z_0} - I\right) \qquad (5)$$

where $z_0$ is a normalized impedance (here, $z_0 = 50$ Ohm) and $I$ denotes an $N \times N$ identity matrix. It should be reiterated that, for the IMPM model, the contents of $x$ only affect the composition of $Z$ in (3), which is obtained as a result of a single multi-port EM simulation. Once $Z$ matrix is extracted, the responses of (5) are only affected by the composition of $y$. It is worth noting that the cost of deriving (3) is higher compared to evaluation of conventional antenna EM models [8]. Computational overhead is associated with the need to extract $Z$ matrix for multiple ports rather than only one (antennas are predominantly single-port devices) [1], [4], [5]. For more comprehensive discussion on IMPM representation of pixel-based structures, see [8], [9].

### 2.3   Optimization Engine for Local Tuning of Pixel Antenna Geometry

Despite low cost, IMPM model of Section 2.2 is suitable only for adjustment of interconnections between the individual pixels. The process can be performed, e.g., using population-based metaheuristics, or exhaustive search [4]-[5]. However, fine-tuning of the topology involves optimization of floating-point parameters $x$ (e.g., dimensions of individual pixels, separation between them, or configuration of the fed line). Conventional algorithms require dozens of EM simulations for convergence, even when $x$ is relatively low-dimensional [5]. Here, the problem is mitigated using a gradient-based method embedded within a trust-region (TR) loop. Let $R_f(x,y^*)$ be the equivalent of a pixel-based EM model which substitutes open/closed ports specified in $y^*$ by metal/etched connections between pixels (cf. Fig. 2(a)). TR framework generates a series of approximations ($i = 0, 1, 2, \ldots$) to the final design $x^*$ by solving [4], [10]:

$$x^{(i+1)} = \arg \min_{\|x - x^{(i)}\| \leq \delta} \left( U\left(R_{sx}^{(i)}(x)\right)\right) \qquad (6)$$

where $R_{sx}^{(i)}(x) = R_f(x^{(i)}, y^*) + J_f(x^{(i)}, y^*)(x - x^{(i)})$ is a first-order Taylor model and $J_f$ is a Jacobian based on large-step finite-differences around $R_f(x^{(i)}, y^*)$ [10]. The trust-region radius $\delta$, i.e., the range around $x^{(i)}$ for which the model $R_{sx}^{(i)}$ is considered acceptably accurate, is controlled based on a gain ratio given as [10]:

$$\rho = \frac{U\left(R_f\left(x^{(i+1)}, y^*\right)\right) - U\left(R_f\left(x^{(i)}, y^*\right)\right)}{U\left(R_{sx}\left(x^{(i+1)}\right)\right) - U\left(R_{sx}\left(x^{(i)}\right)\right)} \qquad (7)$$

At first iteration ($i = 0$), the radius is set to $\delta = 1$. Then, it is adjusted based on standard rules, i.e., $\delta = 2\delta$, when $\rho > 0.75$ and $\delta = \delta/3$, when $\rho < 0.25$. The candidate design obtained from (6) is accepted when $\rho > 0$; otherwise it is rejected. The algorithm is terminated either when $\delta < \varepsilon$, or $\|x^{(i+1)} - x^{(i)}\| < \varepsilon$ (here, $\varepsilon = 10^{-2}$ is used). The computational cost of the algorithm is just $D + 1$ EM simulations per successful iteration ($D$ denotes dimensionality of $x$). Additional EM evaluation is required for each unsuccessful step. Upon termination of the algorithm (6), $x^* = x^{(i)}$ is set. For more comprehensive discussion on TR-based optimization of antennas, see [4], [10], [11].

### 2.4 Feature-Assisted Representation of Antenna Responses

Antenna can be considered as a transformer between the transmission line (e.g., microstrip) and wireless propagation medium. Consequently, operational frequencies of the radiator are proportional to its physical dimensions [11]. The problem outlined in Section 2.1 involves adjustment of $y$ for the given $x$ (that determines size of the radiator) followed by fine tuning of the latter. Due to sequential nature, the procedure might be cumbersome when determination of target operational frequencies is associated with tuning of $x$. The problem is important for multi-band antennas [4]. It can be mitigated by execution of local-search within the domain of features [12]. The latter represent antenna responses in the form of carefully selected points.

Let $F(x) = [\omega\ L] = P(R_f(x))$ be the feature-based response of a pixel-based antenna extracted from EM simulation using a function $P$ [4], [12]. Here, $\omega = [\omega_1 \ldots \omega_q]^T$ and $L = [L_1 \ldots L_q]^T$ ($q = 1, \ldots, Q$) denote the features defined w.r.t. frequencies (typically in GHz) and levels (in dB) of interest for the design problem. The pair of points can represent frequency at the specific response level (e.g., local minimum pertinent to resonance of narrowband antenna), or level at specific target frequency (e.g., at the expected edge of the operational bandwidth) [12]. When compared against $R_f(x)$, changes of response features are a much less non-linear function of $x$. Therefore, coordinate points can be used for reliable shifting of resonances over a wide frequency range while maintaining acceptable cost of the process. Representation of antenna responses using feature points is illustrated in Fig. 2. For more comprehensive discussion on the concept, see [4], [12].

It should be emphasized that, when used with optimization engine of Section 2.3, the $R_{sx}$ model is not only constructed based on $F(x)$ responses—instead of $R_f(x)$—but also evaluated using different design objectives. The same applies to evaluation of gain ratio.

### 2.5 Summary of the Bi-Stage Optimization Method

Unsupervised generation and tuning of the pixel-based topology enables streamlined development of antennas according to the desired design specifications. The considered bi-stage optimization method can be summarized as follows.

1. Set $N_x$, $N_y$, $x$, calculate $M$ and generate antenna EM model;
2. Generate $Z(x, f)$ and construct $R_{sy}(y)$ as explained in Section 2.2;
3. Find $y^*$ by solving (2) using objective $U$ and global-search method of choice;
4. For frequency-based design, go to Step 6; otherwise go to Step 5;
5. Define features extraction function $P$ and substitute $R_f(x)$ with $F(x)$;
6. Find $x^*$ by minimization of (6) using selected objective function and END.

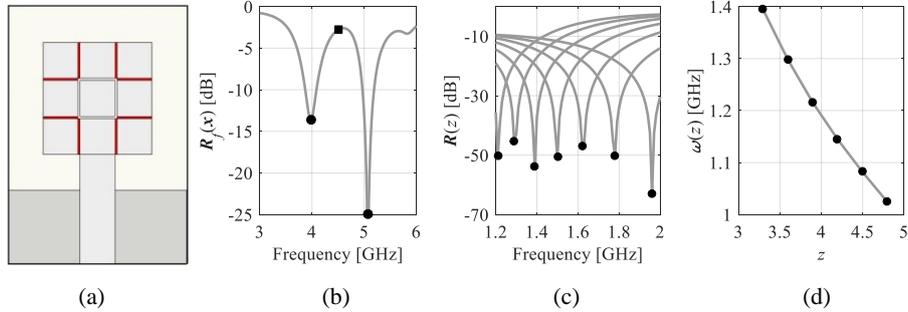

| (a) | (b) | (c) | (d) |

**Fig. 2.** Pixel-based antenna: (a) example EM model with highlight on metal/etched (red/white) connections between pixels, as well as (b) feature-based representation ($Q = 3$) of $R_f(x)$ response (——) in terms of levels (●) and frequencies of interest (■), (c) frequency sweep along given $z$, and (d) feature-based sweep. Changes of $\omega(z)$ are much less non-linear compared to $R(z)$ [4].

Antenna EM model is generated without user inference using a custom script that determines the lattice of pixels constituting the radiator, as well as the ports that enable controlling connections between them. It should be noted that computational cost of the outlined algorithm is low as identification of $x^*$ and $y^*$ involves only a handful of EM simulations.

## 3    Results

The considered framework is demonstrated based on monopole topology constituted by a lattice of pixels with $N_x = 3$ rows and $N_y = 3$ columns, as well as $M = 12$ auxiliary ports (cf. Section 2.2). The EM models incorporate an FR 4 substrate ($h = 1.6$ mm, $\varepsilon_r = 4.3$, $\tan\delta = 0.025$). The initial parameters are $x^{(0)} = [3\ 0.2\ 0\ 3]^T$. Two test cases that involve development of broadband and dual-band radiators are considered. EM models are generated using CST Studio. Given relatively low number of internal ports and negligible cost of $R_{sy}$, minimization of (2) is performed using exhaustive search. Lower and upper bounds on $x$ are given as $l_x = [3\ 0.2\ 0\ 2.4]^T$ and $u_x = [5\ 0.6\ 4\ 5]^T$.

### 3.1    Broadband Antenna Design

The first case concerns development of wideband structure dedicated to work within $f_L = 3.8$ GHz to $f_H = 10$ GHz range. Upon generation of the antenna EM model, the impedance matrix $Z$ is obtained for the design $x^{(0)}$ as explained in Section 2.2. The optimization is governed by a min-max function $U = \max\{\max(R - R_{thd}), 0\}_{fL \leq f \leq fH}$, where $R_{thd} = -10$ dB is selected as a threshold for an acceptable in-band reflection; $f \in f$ (cf. Section 2.1). The pixels configuration $y^* = [\infty\ \infty\ 0\ 0\ 0\ 0\ 0\ \infty\ \infty\ 0\ \infty\ 0]^T$ is found through minimization of (2) using the above objective. Next, the identified geometry is fine-tuned using algorithm of Section 2.3. The final solution $x^* = [3.55\ 0.2\ 1.78\ 2.8]^T$ is found after just 6 iterations of (6), which corresponds to a total of 31 $R_f(x)$ simulations. Overall, the computational cost of antenna development amounts to 33.3 $R_f$ simulations (~0.56 hours of CPU-time on a dual AMD EPYC 7282 system). Figure 3(a) shows a geometry of the optimized antenna and a comparison of its responses at each stage of the design process. It is worth noting that $R_f(x^{(0)},y^*)$ slightly violates the imposed design requirements. However, the performance is improved as a result of the second stage of design process. Overall, $R_f(x^*,y^*)$ is characterized by the reflection below the –10 dB level within 3.74 GHz to 10 GHz range.

### 3.2 Multi-Band Antenna Design

The second case involves design of a dual-band antenna dedicated to work at $f_1 = 3$ GHz and $f_2 = 6$ GHz. Given that $x^{(0)}$ might be inadequate for obtaining resonances at target frequencies, the problem has been reformulated to identify the geometry for which the second frequency is shifted w.r.t. the first resonance by the predefined scaling factor $K$ (here, $K = f_2/f_1 = 2$). For such a problem, identification of suitable topology can be performed by minimization of the following objective function $U = \max\{\mathbf{R} - R_{thd}\}_{fr.1, fr.2 \in f}$, where $f_{r.1}$ is the first resonance frequency featuring the reflection below $R_{thd} = -15$ dB level and $f_{r.2} = Kf_{r.1}$. Similarly as in Section 3.1, the optimum solution $y^* = [0 \infty 0 0 \infty 0 0 \infty \infty 0 0]^T$ is found as a result of exhaustive search. It is worth noting that the search space defined by $U$ is subject to discontinuities. Nonetheless, it does not poses a challenge given the nature of selected search mechanism. As shown in Fig. 3(b), the center frequencies for the design $R_f(x^{(0)}, y^*)$ are at 3.74 GHz and 9.1 GHz, which represents substantial shift w.r.t. the target values. Having that in mind, the structure response is converted to features $F(x)$ similar as the ones shown in Fig. 2(b), i.e., frequencies/levels at resonances. The feature-based objective function is $U_F = \beta_1 \max(\mathbf{L} - R_{tF}) + \|\omega - [f_1 \; f_2]^T\|$, where $\beta_1 = 10$ is used to balance the frequency- and level-related requirements. The design goal involves shifting the resonances to target frequencies, while maintaining the reflection at $\omega$ below the $R_{thd}$ level. The optimized design $x^* = [4.88 \; 0.6 \; 0.79 \; 3.41]^T$ is found after 8 TR iterations—minimization of (6) is realized based on $F(x)$ responses—which corresponds to 33 $R_f(x, y^*)$ evaluations. As shown in Fig. 3(b), the optimized design is characterized by reflection below –10 dB level within 2.66 GHz to 4.06 GHz and 4.44 GHz to 6.67 GHz ranges, respectively. The overall computational cost of the design process corresponds to 35.3 $R_f$ simulations (~0.6 hours of CPU-time) which is low, given identification of $y^*$ based on the exhaustive global-search.

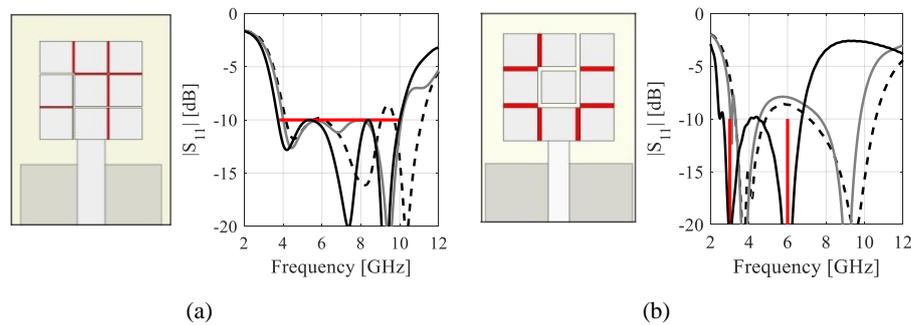

**Fig. 3.** Pixel antennas – optimized geometries and responses for (gray) IMPM $R_{sy}(x^{(0)}, y^*)$, as well as (black) EM models $R_f(x^{(0)}, y^*)$ (– –), and $R_f(x^*,y^*)$ (—): (a) first and (b) second case.

## 4  Conclusion

In this work, a bi-stage framework for cost-efficient development of pixel-based antennas has been considered. The method involves identification of the desirable antenna topology through global optimization of the cost-efficient IMPM model. The resulting geometry is then subject to fine-tuning through adjustment of floating-point dimensions (e.g., size, and separation of pixels). The method has been demonstrated based on two examples concerning automatic design of broadband and dual-band antennas. For the first case, fine-tuning has been executed using frequency responses. For the second case, structure geometry has been adjusted using the feature-based responses to enable shifting the off-target resonances obtained as a result of global-search. For the considered test cases, the aggregated cost of antenna design does not exceed 36 EM simulations.

**Acknowledgments.** This work was supported in part by the National Science Center of Poland Grant 2021/43/B/ST7/01856.

**Disclosure of Interests.** The authors declare no conflicts of interest.